\documentclass[openany, a4paper, 12pt]{amsart}

\usepackage{amsfonts}

 \usepackage{mathrsfs,amsfonts,amsmath,amssymb,amsthm}
 \usepackage{dsfont}
 \usepackage{color}
 \usepackage{braket}
 \usepackage{stmaryrd}

 \makeatletter
 \@addtoreset{equation}{section}
 \makeatother
 \newtheorem{theorem}{Theorem}[section]
 \newtheorem{lemma}[theorem]{Lemma}
 \newtheorem{corol}[theorem]{Corollary}
 \newtheorem{prop}[theorem]{Proposition}

 \def\blemma{\begin{lemma}}\def\elemma{\end{lemma}}
 \def\bproposition{\begin{prop}}\def\eproposition{\end{prop}}
 \def\btheorem{\begin{theorem}}\def\etheorem{\end{theorem}}
 \def\bcorollary{\begin{corol}}\def\ecorollary{\end{corol}}

 \def\beqlb{\begin{eqnarray}}\def\eeqlb{\end{eqnarray}}
 \def\beqnn{\begin{eqnarray*}}\def\eeqnn{\end{eqnarray*}}

 \def\ar{\!\!\!&}

 \def\e{\mbox{\rm e}}

 \def\qqquad{\qquad\qquad}

 \def\proof{\noindent{\it Proof.~~}}\def\qed{\hfill$\Box$\medskip}

 \def\<{\langle}\def\>{\rangle}

 \def\mbb{\mathbb}

 \def\ar{\!\!&}

\begin{document}

\bigskip\bigskip

\centerline{\huge\bf Some limit theorems for flows}

\medskip

\centerline{\huge\bf of branching processes\footnote{ Supported by
NSFC and 985 Project.}}

\bigskip

\centerline{Hui He and Rugang Ma\footnote{Corresponding author. E-mail:
marugang@mail.bnu.edu.cn.}}

\bigskip

\centerline{Beijing Normal University}

\bigskip

{\narrower{\narrower

\noindent\textit{Abstract:} We construct two kinds of stochastic flows of
discrete Galton-Watson branching processes. Some scaling limit theorems
for the flows are proved, which lead to local and nonlocal branching
superprocesses over the positive half line.

\medskip


\medskip


\par}\par}

\bigskip


\section{Introduction}

Continuous-state branching processes (CB-processes) arose as weak limits
of rescaled discrete Galton-Watson branching processes; see, e.g.,
Ji\v{r}ina (1958) and Lamperti (1967). Continuous-state branching
processes with immigration (CBI-processes) are generalizations of them
describing the situation where immigrants may come from other sources of
particles. Those processes can be obtained as the scaling limits of
discrete branching processes with immigration; see, e.g., Kawazu and
Watanabe (1971) and Li (2006). A CBI-process was constructed in Dawson and
Li (2006) as the strong solution of a stochastic equation driven by
Brownian motions and Poisson random measures; see also Fu and Li (2010). A
similar construction was given in Li and Ma (2008) using a stochastic
equation driven by time-space Gaussian white noises and Poisson random
measures.

In the study of scaling limits of coalescent processes with multiple
collisions, Bertoin and Le Gall (2006) constructed a flow of jump-type
CB-processes as the weak solution flow of a system of stochastic equations
driven by Poisson random measures; see also Bertoin and Le Gall (2003,
2005). A more general flow of CBI-processes was constructed in Dawson and
Li (2012) as strong solutions of stochastic equations driven by Gaussian
white noises and Poisson random measures. The flows in Bertoin and Le Gall
(2006) and Dawson and Li (2012) were also treated as path-valued processes
with independent increments. Motivated by the works of Aldous and Pitman
(1998) and Abraham and Delmas (2010) on tree-valued Markov processes,
another flow of CBI-processes was introduced in Li (2012), which was
identified as a path-valued branching process. From the flows in Bertoin
and Le Gall (2006), Dawson and Li (2012) and Li (2012), one can define
some superprocesses or immigration superprocesses over the positive half
line with local and nonlocal branching mechanisms. To study the genealogy
trees for critical branching processes conditioned on non-extinction,
Bakhtin (2011) considered a flow of continuous CBI-processes driven by a
time-space Gaussian white noise. He obtained the flow as a rescaling limit
of systems of discrete Galton-Watson processes and also pointed out the
connection of the model with a superprocess conditioned on non-extinction.

In this paper, we consider two flows of discrete Galton-Watson branching
processes and show suitable rescaled sequences of the flows converge to
the flows of Dawson and Li (2012) and Li (2012), respectively. The main
motivation of the work is to understand the connection between discrete
and continuum tree-valued processes. Our results generalize those of
Bakhtin (2011) to flows of discontinuous CB-processes. To simplify the
presentation, we only treat models without immigration, but the arguments
given here carry over to those with immigration. We shall first prove
limit theorems for the induced superprocesses, from which we derive the
convergence of the finite-dimensional distributions of the path-valued
branching processes.

In Section 2, we give a brief review of the flows of Dawson and Li (2012)
and Li (2012). In Section 3, we consider flows consisting of independent
branching processes and show their scaling limit gives a flow of the type
of Dawson and Li (2012). The formulation and convergence of interactive
flows were discussed in Section 4, which lead to a flow in the class
studied in Li (2012).

Let $\mbb N$ = $\{0,1,2,\cdots\}$ and $\mbb N_+$= $\{1,2,\cdots\}$. For
any $a\geq 0$ let $M[0,a]$ be the set of finite Borel measures on $[0,a]$
endowed with the topology of weak convergence. We identify $M[0,a]$ with
the set $F[0,a]$ of positive right continuous increasing functions on
$[0,a]$. Let $B[0,a]$ be the Banach space of bounded Borel functions on
$[0,a]$ endowed with the supremum norm $\|\cdot\|$. Let $C[0,a]$ denote
its subspace of continuous functions. We use $B[0,a]^+$ and $C[0,a]^+$ to
denote the subclasses of positive elements and $C[0,a]^{++}$ to denote the
subset of $C[0,a]^+$ of functions bounded away from zero. For $\mu\in
M[0,a]$ and $f\in B[0,a]$ write $\<\mu, f\>$ = $\int f d\mu$ if the
integral exists.


\section{Local and nonlocal branching flows}

\setcounter{equation}{0}

In this section, we recall some results on constructions and
characterizations of the flow of CB-processes and the associated
superprocess. It is well-known that the law of a CB-process is determined
by its \emph{branching mechanism} $\phi$, which is a function on
$[0,\infty)$ and has the representation
 \beqlb\label{1.1}
\phi(z)=bz+\frac{1}{2}\sigma^2z^2 +\int_0^\infty (e^{-zu}-1+zu)m(du),
 \eeqlb
where $\sigma\geq 0$ and $b$ are constants and $(u\wedge{u}^2)m(du)$ is a
finite measure on $(0,\infty)$. Let $W(ds,du)$ be a white noise on
$(0,\infty)^2$ based on $dsdu$ and $\tilde{N}(ds,dz,du)$ a compensated
Poisson random measure on $(0,\infty)^3$ with intensity $dsm(dz)du$. By
Theorem~3.1 of Dawson and Li (2012), a CB-process with branching mechanism
$\phi$ can be constructed as the pathwise unique strong solution $\{Y_t:
t\geq 0\}$ to the stochastic equation:
 \beqlb\label{1.02}
Y_t \ar=\ar Y_0 + \sigma\int_0^t\int_0^{Y_{s-}}W(ds,du) - \int_0^t
bY_{s-}ds \cr
 \ar\ar
+ \int_0^t\int_0^{\infty}\int_0^{Y_{s-}}z \tilde{N}(ds,dz,du).
 \eeqlb

Let us fix a constant $a\geq 0$ and a function $\mu\in F[0,a]$. Let
$\{Y_t(q): t\geq 0\}$ denote the solution to (\ref{1.02}) with $Y_0(q) =
\mu(q)$. We can consider the solution flow $\{Y_t(q): t\geq 0, q\in
[0,a]\}$ of (\ref{1.02}). As observed in Dawson and Li (2012), there is a
version of the flow which is increasing in $q\in [0,a]$. Moreover, we can
regard $\{(Y_t(q))_{t\geq 0}: q\in [0,a]\}$ as a path-valued stochastic
process with independent increments. Let $\{Y_t: t\ge 0\}$ denote the
$M[0,a]$-valued process so that $Y_t[0,q] = Y_t(q)$ for every $t\geq 0$
and $q\in [0,a]$. Then $\{Y_t: t\ge 0\}$ is a c\`{a}dl\`{a}ag superprocess with
branching mechanism $\phi$ and trivial spatial motion; see Theorems~3.9
and~3.11 in Dawson and Li (2012). For $\lambda\ge 0$ let $t\mapsto
v(t,\lambda)$ be the unique locally bounded positive solution of
 \beqlb\label{2.4}
v(t,\lambda)= \lambda-\int_0^t\phi(v(s,\lambda))ds, \qquad t\geq 0.
 \eeqlb
For any $f\in B[0,a]^+$ define $x\mapsto v(t,f)(x)$ by $v(t,f)(x) =
v(t,f(x))$. Then the superprocess $\{Y_t: t\ge 0\}$ has transition
semigroup $(Q_t)_{t\geq 0}$ on $M[0,a]$ defined by
 \beqlb\label{2.5}
\int_{M[0,a]} e^{-\<\nu,f\>}Q_t(\mu,\nu)
 =
\exp\left\{-\<\mu,v(t,f)\>\right\}, \qquad f\in B[0,a]^+.
 \eeqlb
By Proposition~3.1 in Li (2011) one can see that $v(t,f)\in C[0,a]^{++}$
for every $f\in C[0,a]^{++}$. Then it is easy to verify that $(Q_t)_{t\geq
0}$ is a Feller semigroup.

We can define another branching flow. For this purpose, let us consider an
admissible family of branching mechanisms $\{\phi_q: q\in [0,a]\}$, where
$\phi_q$ is given by (\ref{1.1}) with parameters $(b,m)=(b_q,m_q)$
depending on $q\in [0,a]$. Here by an \emph{admissible family} we mean for
each $z\geq 0$, the function $q\mapsto \phi_q(z)$ is decreasing and
continuously differentiable with the derivative $\psi_\theta(z) =
-(\partial/\partial\theta)\phi_{\theta}(z)$ of the form
 \beqlb\label{1.2}
\psi_\theta(z)=h_{\theta}z+\int_0^{\infty}(1-e^{-zu})n_{\theta}(du),
 \eeqlb
where $h_{\theta}\geq 0$ and $n_{\theta}(du)$ is a $\sigma$-finite kernel
from $[0,a]$ to $(0,\infty)$ satisfying
 \beqnn
\sup_{0\leq\theta\leq
a}\Big[h_{\theta}+\int_0^{\infty}un_{\theta}(du)\Big]<\infty.
 \eeqnn
Then we have
 \beqlb\label{1.3}
\phi_q(z)=\phi_0(z)-\int_0^q \psi_{\theta}(z)d\theta, \qquad z\geq 0.
 \eeqlb
Let $m(dz,d\theta)$ be the measure on $(0,\infty)\times[0,a]$ defined by
 \beqnn
m([c,d]\times [0,q])=m_q[c,d],\quad q\in [0,a], d>c>0.
 \eeqnn
Suppose that $W(ds,du)$ is a white noise on $(0,\infty)^2$ based on $dsdu$
and $\tilde{N}(ds,dz,d\theta,du)$ is a compensated Poisson random measure
on $(0,\infty)^2\times [0,a]\times (0,\infty)$ with intensity
$dsm(dz,d\theta)du$. By the results in Li (2012), for any $\mu\in F[0,a]$
the stochastic equation
 \beqlb\label{4.16}
Y_t(q) \ar=\ar \mu(q)-b_q\int_0^t Y_{s-}(q)ds+
\sigma\int_0^t\int_0^{Y_{s-}(q)}W(ds,du)\cr
 \ar\ar
 +\int_0^t\int_0^{\infty}\int_{[0,q]}\int_0^{Y_{s-}(q)}
 z \tilde{N}(ds,dz,d\theta,du)
 \eeqlb
has a unique solution flow $\{Y_t(q): t\geq 0,q\in [0,a]\}$. For each
$q\in [0,a]$, the one-dimensional process $\{Y_t(q): t\geq 0\}$ is a
CB-process with branching mechanism $\phi_q$. It was proved in Li (2011)
that there is a version of the flow which is increasing in $q\in [0,a]$.
Moreover, we can also regard $\{(Y_t(q))_{t\geq 0}: q\in [0,a]\}$ as a
path-valued branching process. The solution flow of (\ref{4.16}) also
induces a c\`{a}dl\`{a}ag superprocess $\{Y_t: t\ge 0\}$ with state space
$M[0,a]$. Let $f\mapsto \Psi(\cdot,f)$ be the operator on $C^+[0,a]$
defined by
 \beqlb\label{1.4}
\Psi(x,f)=\int_{[0,a]} f(x\vee \theta)h_{\theta}d \theta
+\int_{[0,a]}d\theta\int_0^{\infty}(1-e^{-zf(x\vee
\theta)})n_{\theta}(dz).
 \eeqlb
The superprocess $\{Y_t: t\ge 0\}$ has local branching mechanism $\phi_0$
and nonlocal branching mechanism given by (\ref{1.4}); see Theorem~6.2 in
Li (2012). Then the transition semigroup $(Q_t)_{t\geq 0}$ of $\{Y_t: t\ge
0\}$ is defined by
 \beqlb\label{1.5}
\int_{M[0,a]}e^{-\<\nu,f\>}Q_t(\mu,d\nu)=\exp\Big\{-\<\mu,V_tf\>\Big\},
\qquad f\in C^+[0,a],
 \eeqlb
where $t\mapsto V_tf$ is the unique locally bounded positive solution of
 \beqlb\label{1.6}
V_tf(x)=f(x)-\int_0^t [\phi_0(V_sf(x))-\Psi(x,V_sf)]ds,\qquad t\geq 0,
x\in [0,a].
 \eeqlb

To study the scaling limit theorems of the discrete branching flows, we
need to introduce a metric on $M[0,a]$. Let $\{h_0,h_1,h_2,\cdots\}$ be a
countable dense subset of $\{h\in C[0,a]^+:\|h\|\leq 1\}$ with $h_0\equiv
1$. For convenience we assume each $h_i$ is bounded away from zero. Then
$\{h_0,h_1,h_2,\cdots\}\subset C[0,a]^{++}$. Now we define a metric $\rho$
on $M[0,a]$ by
 \beqnn
\rho (\mu,\nu) = \sum_{i=0}^{\infty}\frac{1}{2^i}(1\wedge|\<\mu,h_i\>
-\<\nu,h_i\>|), \qquad \mu, \nu\in M[0,a].
 \eeqnn
It is easy to see that the metric is compatible with the weak convergence
topology of $M[0,a]$. In other words, we have $\mu_n\to\mu$ in $M[0,a]$ if
and only if $\rho(\mu_n,\mu)\to 0$. For $\nu\in M[0,a]$, set
$\e_{h_i}(\nu)= \e^{-\<\nu,h_i\>}$.

\btheorem\label{t2.1} The metric space $(M[0,a], \rho)$ is a locally
compact Polish (complete and separable) space, and
$\{\e_{h_{i}}:i=0,1,2,\cdots\}$ is strongly separating points of $M[0,a]$,
that is, for every $\nu\in M[0,a]$ and $\delta>0$, there exists a finite
set $\{\e_{h_{i_1}},\e_{h_{i_2}},\cdots,\e_{h_{i_k}}\}\subset
\{\e_{h_i}:i=0,1,2,\cdots\}$ such that
 \beqnn
\inf_{\mu:\rho(\mu,\nu)\geq \delta}\max_{1\leq j\leq
k}|\e_{h_{i_j}}(\mu)-\e_{h_{i_j}}(\nu)|>0.
 \eeqnn
 \etheorem

\proof By Li (2011, p.4 and p.7) we know $M[0,a]$ is separable and locally
compact, so there is a complete metric on $M[0,a]$ compatible with the
weak convergence topology. The following argument shows the metric $\rho$
defined above is complete. Suppose $\{\mu_n\}_{n\geq 1}\subset M[0,a]$ is
a Cauchy sequence under $\rho$. Then for every $m\geq 1$,
$\{\<\mu_n,h_m\>\}_{n\geq 1}$ is also a Cauchy sequence. We denote the
limit by $\Phi(h_m)$. For $f\in C[0,a]^+$ satisfying $\|f\|\leq 1$, let
$\{h_{i_k}\}_{k\geq 1}\subset \{h_0,h_1,h_2,\cdots\}$ be a sequence so
that $\|h_{i_k}-f\|\to 0$ as $k\to\infty$. For $n\geq m\geq 1$ we have
 \beqnn
\limsup_{m,n\to\infty}|\<\nu_n,f\>-\<\nu_m,f\>|
 \ar\leq\ar
\limsup_{m,n\to\infty}\Big[|\<\nu_n,f\>-\<\nu_n,h_{i_k}\>|\cr
 \ar\ar
+|\<\nu_n,h_{i_k}\>-\<\nu_m,h_{i_k}\>|+|\<\nu_m,h_{i_k}\>-\<\nu_m,f\>|\Big]\cr
 \ar\leq\ar
 2\Phi(1)\|f-h_{i_k}\|.
 \eeqnn
Then letting $k\to\infty$ we have
 \beqnn
\limsup_{m,n\to\infty}|\<\nu_n,f\>-\<\nu_m,f\>|=0.
 \eeqnn
By linearity the above relation holds for all $f\in C[0,a]$, so the limit
$\Phi(f)=\lim\limits_{n\to\infty}\<\mu_n,f\>$ exists for each $f\in
C[0,a]$. Clearly, $f\to\Phi(f)$ is a positive linear functional on
$C[0,a]$. By the Riesz representation theorem there exists $\mu\in M[0,a]$
so that $\<\mu,f\> =\Phi(f)$ for every $f\in C[0,a]$. By the construction
of $\Phi$ we have $\mu_n\to\mu$, therefore, $\rho(\mu_n,\mu)\to 0$. That
proves the first assertion of the theorem.

For any $\nu\in M[0,a]$ and $\delta\geq 0$, there exists an $N_0\in\mbb
N_+$ such that $\sum_{i=N_0+1}^{\infty}1/{2^i}<{\delta}/{2}$. Consider
$\{h_0,h_1,\cdots,h_{N_0}\}$, for any $\mu\in M[0,a]$ satisfying
{$\rho(\mu,\nu)\geq \delta$}, we have
 \beqnn
\sum_{i=0}^{N_0}\frac{1}{2^i}(1\wedge|\<\mu,h_i\>-\<\nu,h_i\>|)\geq
\frac{\delta}{2},
 \eeqnn
and thus,
 \beqnn
\sum_{i=0}^{N_0}(1\wedge|\<\mu,h_i\>-\<\nu,h_i\>|)\geq \frac{\delta}{2}.
 \eeqnn
It follows that
 \beqnn
|\<\mu,h_j\>-\<\nu,h_j\>|\geq \frac{\delta}{2N_0}~~\mbox{for some}~~0\leq
j\leq N_0.
 \eeqnn
Since
 \beqnn
|\e^{-x}-\e^{-y}|=\e^{-y}|\e^{y-x}-1|\geq \e^{-y}\Big[(\e^{|y-x|}-1)\wedge
(1-\e^{-|y-x|})\Big],\qquad x,y\in \mbb R,
 \eeqnn
 we have
 \beqnn
\ar\ar\inf_{\mu:\rho(\mu,\nu)\geq \delta}\max_{0\leq i\leq
N_0}|\e_{h_{i}}(\mu)-\e_{h_{i}}(\nu)|\cr
 \ar\ar\qquad\quad
\geq\e^{-\max_{0\leq i\leq
N_0}\<\nu,h_i\>}\Big[(\e^{\frac{\delta}{2N_0}}-1)\wedge
(1-\e^{-\frac{\delta}{2N_0}})\Big]\cr
 \ar\ar\qquad\quad
> 0.
 \eeqnn
That proves the second assertion. \qed


\section{Flows of independent branching processes}

\setcounter{equation}{0}

In this section, we consider some flows of independent Galton-Watson
branching processes. We shall study the scaling limit in the setting of
superprocesses. Then we derive the convergence of the finite-dimensional
distributions of the path-valued processes.

Let $\{g_i: i=0,1,2,\cdots\}$ be a family of probability generating
functions. Given a family of $\mbb{N}$-valued independent random variables
$\{X_0(i): i=0,1,2,\cdots\}$, for each $i\in\mbb N$ suppose that there are
$X_0(i)$ independent Galton-Watson trees originating at time $0$ and at
place $i$ with offspring distribution given by $g_i$. Let us denote by
$X_n(i)$ the numbers of vertices in the $n$-th generation of the trees
with root at $i$. In addition, we assume $(X_n(i))_{n\geq 0}$,
$i=1,2,\cdots$ are mutually independent. It is well-known that for each
$i\in\mbb N$, $(X_n(i))_{n\geq 0}$ is a \emph{Galton-Watson
branching process} (GW-process) with parameter $g_i$; i.e., a discrete-time $\mbb{N}$-valued
Markov chain with $n$-step transition matrix $P^n(j,k)$ defined by
 \beqlb\label{2.1}
 \sum_{k=0}^{\infty}P^n(j,k)z^k=(g_i^n(z))^j,\qquad |z|\leq 1,
 \eeqlb
where $g_i^n(z)$ is defined by $g_i^n(z)=g_i(g_i^{n-1}(z))$ successively
with $g_i^0(z)=z$.

Suppose that for each integer $k\geq 1$ we have a sequence of GW-processes
$\{(X_n^{(k)}(i))_{n\geq 0}:i\geq 0\}$ with parameter $g_i^{(k)}$. Let
$\gamma_k$ be a positive real sequence so that $\gamma_k\to\infty$
increasingly as $k\to\infty$. For $m,n\in \mbb N$, define
 \beqnn
\bar{X}_n^{(k)}(m)=\sum_{i=0}^m X_n^{(k)}(i),
 \eeqnn
and
 \beqnn
\qqquad\qqquad
Y_t^{(k)}(x)=\displaystyle\frac{1}{k}\bar{X}_{[\gamma_kt]}^{(k)}([kx]),
\qquad k=1,2,\cdots,
 \eeqnn
where $[\cdot]$ denotes the integer part. Then the increasing function
$x\mapsto Y_t^{(k)}(x)$ induces a random measure $Y_t^{(k)}(dx)$ on
$[0,\infty)$ so that $Y_t^{(k)}([0,x])=Y_t^{(k)}(x)$ for $x\geq 0$. For
convenience we fix a constant $a\geq 0$ and consider the restriction of
$\{Y_t^{(k)}:t\geq 0\}$ to $[0,a]$ without changing the notation. Clearly,
 \beqnn
Y_0^{(k)}=\frac{1}{k}\sum_{i=0}^{[ka]}X_0^{(k)}(i)\delta_{\frac{i}{k}}
 \eeqnn
and
 \beqnn
Y_t^{(k)}=\frac{1}{k}\sum_{i=0}^{[ka]}X_{[\gamma_kt]}^{(k)}(i)\delta_{\frac{i}{k}}.
 \eeqnn
In view of (\ref{2.1}), for each $i\geq 0$, given $X_0^{(k)}(i)=x_i\in
\mbb N$, the conditional distribution $Q_{i,k}^{[\gamma_kt]}(x_i/k,\cdot)$
of $\{k^{-1}X_{[\gamma_kt]}^{(k)}(i):t\geq 0\}$ on
$E_k=\{0,1/k,2/k,\cdots\}$ is determined by
 \beqlb\label{2.2}
\int_{E_k} e^{-\lambda y}Q_{i,k}^{[\gamma_kt]}({x_i}/{k},dy) =
\exp\left\{-\frac{x_i}{k}v_i^{(k)}(t,\lambda)\right\},
 \eeqlb
where $v_i^{(k)}(t,\lambda)=-k\log
(g_i^{(k)})^{[\gamma_kt]}(e^{-\lambda/k})$.

Let $Q_{\mu_k}^{(k)}$ denote the conditional law given
$Y_0^{(k)}=\mu_k=k^{-1}\sum_{i=0}^{[ka]}x_i\delta_{i/{k}}\in M_k[0,a]$,
where $M_k[0,a]:=\{k^{-1}\sum_{i=0}^{[ka]}x_i\delta_{i/{k}}: ~x_i\in{\mbb
N}, ~{k}^{-1}\sum_{i=0}^{[ka]}{x_i}<\infty\}$. For $f\in B[0,a]^+$, from
(\ref{2.2}) we have
 \beqlb\label{2.3}
Q_{\mu_k}^{(k)}\exp\Big\{-\<Y_t^{(k)},f\>\Big\}
 \ar=\ar
Q_{\mu_k}^{(k)}
\exp\bigg\{-\sum_{i=0}^{[ka]}\frac{1}{k}X_{[\gamma_kt]}^{(k)}(i)
f{({i}/{k})}\bigg\}\cr
 \ar=\ar
\prod_{i=1}^{[ka]}\int_{E_k}
e^{-f{({i}/{k})}y}Q_{i,k}^{[\gamma_kt]}({x_i}/{k},dy)\cr
 \ar=\ar
\exp\bigg\{-\sum_{i=0}^{[ka]}\frac{x_i}{k}v_i^{(k)}(t,f({i}/{k}))\bigg\}\cr
 \ar=\ar
\exp\Big\{-\<\mu_k,v^{(k)}(t,f)\>\Big\},
 \eeqlb
where $x\mapsto v^{(k)}(t,f)(x)$ is defined by
$v^{(k)}(t,f)(x)=v_{[kx]}^{(k)}(t,f(x))$.

For any $x, z\geq 0$ define
 \beqlb\label{3.3a}
 \phi_k(x,z)=k\gamma_k[g_{[kx]}^{(k)}(e^{-z/k})-e^{-z/k}].
 \eeqlb
For convenience of statement of the results, we formulate the following
condition:

\noindent{\bf Condition (3.A)} ~For each $a\geq 0$ the sequence
$\{\phi_k(x,z)\}$ is Lipschitz with respect to $z$ uniformly on
$[0,\infty)\times [0,a]$ and there is a continuous function $(x,z)\mapsto
\phi (x, z)$ such that $\phi_k(x, z) \to\phi (x, z)$ uniformly on
$[0,\infty)\times [0,a]$ as $k\to \infty$.

Before giving the limit theorem for the sequence of the rescaled
processes, we first introduce the limit process. By Proposition~4.3 in Li
(2011), if Condition (3.A) is satisfied, the limit function $\phi$ has the
representation
 \beqlb\label{1.01}
\phi(x,z)=b(x)z+\frac{1}{2}c(x)z^2 +\int_0^\infty
(e^{-zu}-1+zu)m(x,du),\qquad x,z\geq 0.
 \eeqlb
where $b$ is a bounded function on $[0,\infty)$ and $c$ is a positive
bounded function on $[0,\infty)$. $(u\wedge{u}^2)m(x,du)$ is a bounded
kernel from $[0,\infty)$ to $(0,\infty)$. Conversely, for any continuous
function $(x,z)\mapsto \phi (x, z)$ given by (\ref{1.01}), we can
construct a family of probability generating functions
$\{g_i^{(k)}:i=0,1,2,\cdots\}$ so that the sequence (\ref{3.3a}) satisfies
Condition (3.A); see, e.g., Li (2011, p.93).

 For any $l\geq 0$, let $B_l[0,\infty)^+$ be the set of positive bounded
functions on $[0,\infty)$ satisfying $\|f\|\leq l$. By a modification of
the proof of Theorem 3.42 in Li (2011), it is not hard to show that for
each $T\geq0$ and $l\geq0$, $v^{(k)}(t,f)(x)$ converges uniformly on the set
$[0,T]\times [0,\infty)\times B_l[0,\infty)^+$ of $(t,x,f)$ to the unique
locally bounded positive solution $(t,x)\mapsto v(t,f)(x)$ of the
evolution equation
 \beqlb\label{1.03}
 v(t,f)(x)=f(x)-\int_0^t
\phi(x,v(s,f)(x))ds. \eeqlb
Let $\{Y_t: t\ge 0\}$ be the superprocess with state space $M[0,a]$
and transition semigroup $(Q_t)_{t\geq 0}$ defined by
 \beqlb\label{1.04}
\int_{M[0,a]} e^{-\<\nu,f\>}Q_t(\mu,\nu)
 =
\exp\left\{-\<\mu,v(t,f)\>\right\}, \qquad f\in B[0,a]^+.
 \eeqlb
Using (\ref{1.03}) and Gronwall's inequality one can see $x\mapsto
v(t,f)(x)$ is continuous on $[0,a]$ for every $f\in C[0,a]^+$. Then by
Proposition~3.1 in Li (2011) it is easy to see that $ v(t,f)\in
C[0,a]^{++}$ for every $f\in C[0,a]^{++}$. From this and (\ref{1.04}) it
follows that $(Q_t)_{t\geq 0}$ is a Feller semigroup. Note that if $\phi
(x, z)= \phi (z)$ independent of $x\geq 0$, then $(Q_t)_{t\geq 0}$ is the
same transition semigroup as that defined by (\ref{2.4}) and (\ref{2.5}).
In this case, the corresponding superprocess can be defined by the
stochastic integral equation (\ref{1.02}).

Let $D([0,\infty),M[0,a])$ denote the space of c\`{a}dl\`{a}g paths from
$[0,\infty)$ to $M[0,a]$ furnished with the Skorokhod topology. The proof
of the next theorem is a modification of that of Theorem~3.43 in Li
(2011).

\btheorem\label{t2.2} Suppose that Condition (3.A) is satisfied. Let
$\{Y_t:t\geq 0\}$ be a c\`{a}dl\`{a}g superprocess with transition
semigroup $(Q_t)_{t\geq 0}$ defined by (\ref{1.03}) and (\ref{1.04}). If
$Y_0^{(k)}$ converges to $Y_0$ in distribution on $M[0,a]$, then
$\{Y_t^{(k)}:t\geq 0\}$ converges to $\{Y_t:t\geq 0\}$ in distribution on
$D([0,\infty),M[0,a])$.\etheorem

\proof For $f\in C[0,a]^{++}$ and $\nu\in M[0,a]$ set
$e_f(\nu)=e^{-\<\nu,f\>}$. Clearly, the function $\nu\mapsto e_f(\nu)$ is
continuous in $\rho$. We denote by $D_1$ the linear span of $\{\e_f:f\in
C[0,a]^{++}\}$. By Theorem \ref{t2.1} we have $D_1$ is an algebra strong
separating the points of $M[0,a]$. Let $C_0(M[0,a])$ be the space of
continuous functions on $M[0,a]$ vanishing at infinity. Then $D_1$ is
uniformly dense in $C_0(M[0,a])$ by the Stone-Weierstrass theorem; see,
e.g., Hewitt and Stromberg (1975, pp.98-99). On the other hand, for any
$f\in C[0,a]^{++}$, since $v(t,f)$ is bounded away from zero and
$v_k(t,f)(x)\to v(t,f)(x)$ uniformly on $[0,\infty)$ for every $t\geq 0$,
we have $v_k(t,f)$ is also bounded away from zero for $k$ sufficiently
large. Without loss of generality we may assume $v_k(t,f)\geq c$ and $v(t,
f)\geq c$ for some $c>0$. Let $Q_t^{(k)}$ denote the transition semigroup
of $Y_t^{(k)}$. We get from (\ref{2.3}) and (\ref{2.5}) that, for any
$M\geq 0$,
 \beqnn
\ar\ar\sup_{\nu\in M_k[0,a]}\left|Q_t^{(k)} e_f(\nu)-Q_t
e_f(\nu)\right|\cr
 \ar\ar\qquad\quad
=\sup_{\nu\in M_k[0,a]}\Big|\exp\Big\{-\<\nu,v_k(t,f)\>\Big\} -
\exp\Big\{-\<\nu,v(t,f)\>\Big\}\Big|\cr
 \ar\ar\qquad\quad
\leq\sup_{\<\nu,1\>\leq M\atop\nu\in M_k[0,a]}
\Big|\exp\Big\{-\<\nu,v_k(t,f)\>\Big\}-\exp\Big\{-\<\nu,v(t,f)\>\Big\}\Big|\cr
\ar\ar\qquad\qquad +\sup_{\<\nu,1\> > M\atop\nu\in M_k[0,a]}
\Big|\exp\Big\{-\<\nu,v_k(t,f)\>\Big\}-\exp\Big\{-\<\nu,v(t,f)\>\Big\}\Big|\cr
 \ar\ar\qquad\quad
\leq\sup_{\<\nu,1\>\leq M\atop\nu\in M_k[0,a]}
|\<\nu,v_k(t,f)\>-\<\nu,v(t,f)\>| + \sup_{\<\nu,1\> > M\atop\nu\in
M_k[0,a]}2e^{-\<\nu,c\>}\cr
 \ar\ar\qquad\quad
\leq M\|v_k(t,f)-v(t,f)\|+2e^{-Mc}.
 \eeqnn
Since $M\geq 0$ was arbitrary, we have
 \beqnn
\lim_{k\to\infty}\sup_{\nu\in M_k[0,a]}\left|Q_t^{(k)} e_f(\nu)-Q_t
e_f(\nu)\right|=0
 \eeqnn
for every $t\geq 0$. Thus
 \beqnn
\lim_{k\to\infty}\sup_{\nu\in M_k[0,a]}\left|Q_t^{(k)} F(\nu)-Q_t
F(\nu)\right|=0
 \eeqnn
for every $t\geq 0$ and $F\in C_0(M[0,a])$. By Ethier and Kurtz (1986,
p.226 and pp.233-234) we conclude that $\{Y_t^{(k)}:t\geq 0\}$ converges
to $\{Y_t:t\geq 0\}$ in distribution on $D([0,\infty),M[0,a])$. \qed

 Let
$\{0\leq a_1<a_2<\cdots<a_n=a\}$ be an ordered set of constants. Denote by
$\{Y_{t,a_i}:t\geq 0\}$ and $\{Y_{t,a_i}^{(k)}:t\geq 0\}$ the restriction of $\{Y_t:t\geq 0\}$ and $\{Y_t^{(k)}:t\geq 0\}$ to $[0,a_i]$, $i=1,2,\cdots,n$, respectively. The following theorem is an extension of Theorem \ref{t2.2}.

\btheorem\label{t2.3} Suppose that Condition (3.A) is satisfied.  If $Y_{0,a}^{(k)}$ converges to $Y_{0,a}$ in
distribution on $M[0,a]$, then $\{(Y_{t,a_1}^{(k)}, \cdots,
Y_{t,a_n}^{(k)}):t\geq 0\}$ converges to $\{(Y_{t,a_1}, \cdots,
Y_{t,a_n}):t\geq 0\}$ in distribution on $D([0,\infty),M[0,a_1]\times
\cdots\times M[0,a_n])$. \etheorem

\proof Let $f_i\in C[0,a_i]$ for $i=1, \cdots, n$. By Theorem~\ref{t2.2}
we see that for every $1\leq i\leq n$, $\{\<Y_{t,a_i}^{(k)},f_i\>: t\geq
0\}$ is tight in $D([0,\infty), \mbb{R})$. Thus
$\{\sum_{i=1}^n\<Y_{t,a_i}^{(k)},f_i\>: t\geq 0\}$ is tight in
$D([0,\infty), \mbb{R})$. Then the tightness criterion of Roelly (1986)
implies $\{(Y_{t,a_1}^{(k)}, \cdots, Y_{t,a_n}^{(k)}):t\geq 0\}$ is tight
in $D([0,\infty),M[0,a_1]\times \cdots\times M[0,a_n])$. Let
$\{(Z_{t,a_1}, \cdots, Z_{t,a_n}):t\geq 0\}$ be a weak limit point of
$\{(Y_{t,a_1}^{(k)}, \cdots, Y_{t,a_n}^{(k)}):t\geq 0\}$. By an argument
similar to the proof of Theorem 5.8 in Dawson and Li (2012) one can show
that $\{(Z_{t,a_1}, \cdots, Z_{t,a_n}):t\geq 0\}$ and $\{(Y_{t,a_1},
\cdots, Y_{t,a_n}):t\geq 0\}$ have the same distributions on
$D([0,\infty),M[0,a_1]\times\cdots\times M[0,a_n])$. That gives the
desired result.\qed

\bcorollary\label{c2.4} Suppose that Condition (3.A) is satisfied. Let
$\{0\leq a_1<a_2<\cdots<a_n=a\}$ be an ordered set of constants. Let
$Y_{t}(a_i):=Y_t[0,a_i]$ and  $Y_{t}^{(k)}(a_i):=Y_t^{(k)}[0,a_i] $ for every $t\geq 0$, $i=1,2,\cdots,n$. If
$(Y_{0}^{(k)}(a_1), \cdots, Y_{0}^{(k)}(a_n))$ converges to $(Y_{0}(a_1),
\cdots, Y_{0}(a_n))$ in distribution on $\mbb{R}_+^{n}$, then
$\{(Y_{t}^{(k)}(a_1), \cdots, Y_{t}^{(k)}(a_n)):t\geq 0\}$ converges to
$\{(Y_{t}(a_1), \cdots, Y_{t}(a_n)):t\geq 0\}$ in distribution on
$D([0,\infty), \mbb{R}_+^{n})$. \ecorollary


\section{Flows of interactive branching processes}

\setcounter{equation}{0}

In this section, we prove some limit theorems for a sequence of flows of
interactive branching processes, which leads to a superprocesses with
local branching and nonlocal branching. From those limit theorems we
derive the convergence of the finite-dimensional distributions of the
path-valued branching processes.

Let $g_0$ be a probability generating function and $\{h_i:i=1,2,\cdots\}$
a family of probability generating functions. For each $i\geq 1$ define
$g_i:=g_0h_1\cdots h_i$ and suppose that $\{\xi_{n,j}(i):
n=0,1,2,\cdots;j=1,2,\cdots\}$ and $\{\eta_{n,j}(i):
n=0,1,2,\cdots;j=1,2,\cdots\}$ are two independent families of positive
integer-valued i.i.d. random variables with distributions given by $g_i$
and $h_i$, respectively. Given another family of positive integer-valued
random variables $\{z_i:i=1,2,\cdots\}$ independent of
$\{\xi_{n,j}(i):i=1,2,\cdots\}$ and $\{\eta_{n,j}(i):i=1,2,\cdots\}$, we
define inductively $X_0(0)=z_0$ and
 \beqlb\label{3.1}
 X_{n+1}(0)=\sum_{j=1}^{X_{n}(0)}\xi_{n,j}(0),
\qquad n=0,1,2,\cdots.
 \eeqlb
Suppose that $\{X_{n}(i):n=0,1,2,\cdots\}$ has been constructed for
$i=0,1,\cdots,m-1$, we define $\{X_{n}(m):n=0,1,2,\cdots\}$ by
$X_0(m)=z_{m}$ and
 \beqlb\label{3.2}
 X_{n+1}(m)=\sum_{j=1}^{X_{n}(m)}\xi_{n,j}(m)
+\sum_{j=1}^{\bar{X}_{n}(m-1)}\eta_{n,j}(m), \qquad n=0,1,2,\cdots,
 \eeqlb
where $\bar{X}_n(m-1)=\sum_{i=0}^{m-1} X_{n}(i),~n=0,1,2,\cdots$. It is
easy to show that for any $m\in {\mbb N}$,
$\{(X_{n}(0),X_{n}(1),\cdots,X_{n}(m)): n=0,1,2,\cdots\}$ is a
discrete-time $\mbb{N}^{m+1}$-valued Markov chain with one-step transition
probability $Q(x,dy)$ determined by, for $\lambda, x\in \mbb{N}^{m+1}$,
 \beqlb\label{3.3}
\int_{\mbb{N}^{m+1}}e^{-\<\lambda, y\>}Q(x,dy) =\prod_{i=0}^m
[g_i(e^{-\lambda_i})]^{x_i}[h_i(e^{-\lambda_i})]^{\sum_{j=0}^{i-1}x_j},
 \eeqlb
where $x_i$ and $\lambda_i$ denote the $i$-th component of $x$ and
$\lambda$, respectively.

Suppose that for each integer $k\geq 1$ we have two sequence of processes
$\{(X_n^{(k)}(i))_{n\geq 0}: {i\geq 0}\}$ and
$\{(\bar{X}_n^{(k)}(i))_{n\geq 0}: i\geq 0\}$ with parameters $g_0^{(k)}$
and $\{h_i^{(k)}: i=1,2,\cdots \}$. Suppose that $\gamma_k$ is a positive
real sequence so that $\gamma_k\to\infty$ increasingly as $k\to\infty$.
Let $[\gamma_kt]$ denote the integer part of $\gamma_kt\geq 0$. Define
 \beqlb\label{3.4}
Y_t^{(k)}(x) :=\frac{1}{k}\bar{X}_{[\gamma_kt]}^{(k)}([kx])
=\frac{1}{k}\sum_{i=0}^{[kx]}{X}_{[\gamma_kt]}^{(k)}(i), \qquad
k=1,2,\cdots.
 \eeqlb
Let $Y_t^{(k)}(dx)$ denote the random measure on $[0,\infty)$ induced by
the random function $Y_t^{(k)}(x)$. We are interested in the asymptotic
behavior of the continuous-time process $\{Y_{t}^{(k)}(dx):t\geq 0\}$ as
$k\to\infty$. Let $h_0^{(k)}\equiv 1$. For any $z\geq 0$ and $\theta\geq
0$ set
 \beqlb\label{4.4a}
\phi_{\theta}^{(k)}(z) = k\gamma_k\Big[g_{[k\theta]}^{(k)}(e^{-z/k})
-e^{-z/k}\Big]
 \eeqlb
and
 \beqlb\label{4.4b}
\psi_{\theta}^{(k)}(z) = k^2\gamma_k\Big[1-h_{[k\theta]}^{(k)}
(e^{-z/k})\Big].
 \eeqlb
Let us consider the following set of conditions:

\noindent\emph{{\bf Condition (4.A)}  ~For every $l\geq 0$, the sequence $\{\phi_0^{(k)}\}$ is
uniformly Lipschitz on $[0,l]$ and there is a function
$\phi_0$ on $[0,\infty)$ such that $\phi_0^{(k)}(z)\to \phi_0 (z)$
uniformly on $[0,l]$ as $k\to \infty$.}

\noindent\emph{{\bf Condition (4.B)} ~ There is a function $\psi$ on
$[0,\infty)^2$ such that, for every $l\geq 0$, $\psi_{\theta}^{(k)}(z)\to
\psi_{\theta}(z)$ uniformly on $[0,l]^2$ as $k\to \infty$ and
 $$
\sup_{\theta\in[0,a]}\frac{d}{dz}
\psi_{\theta}(z)|_{z=0^+}< \infty.
 $$}

\bproposition\label{p3.1} If Conditions (4.A) and (4.B) hold, then for
every $q\geq 0$ there is a branching mechanism $\phi_q$ such that
$\phi_q^{(k)}(z)\to \phi_q (z)$ uniformly on $[0,l]$ for every $l\geq 0$
as $k\to \infty$. Moreover, the family of branching mechanisms $\{\phi_q:
q\geq 0\}$ is admissible with $(\partial/\partial\theta) \phi_{\theta}(z)
= - \psi_{\theta}(z)$. \eproposition

\proof If Conditions (4.A) and (4.B) hold, then the limit function
$\phi_0$ has the representation (\ref{1.1}) with $(b,m)=(b_0,m_0)$ and
$\psi_{\theta}$ has the representation (\ref{1.2}); see, e.g., Li (2011,
p.76). By the definition of $g_i^{(k)}$ it is simple to check that, for
every $q\geq 0$,
 \beqlb\label{3.5}
\phi_q^{(k)}(z)
 \ar=\ar
k\gamma_k[g_0^{(k)}(e^{-z/k})-e^{-z/k}]\prod_{i=1}^{[kq]}
h_i^{(k)}(e^{-z/k})\cr
 \ar\ar
-\sum_{i=1}^{[kq]}k\gamma_k[1-h_i^{(k)}(e^{-z/k})]e^{-z/k}
\prod_{j=i+1}^{[kq]}h_j^{(k)}(e^{-z/k}).
 \eeqlb
By elementary calculations,
 \beqnn
\prod_{i=1}^{[kq]}h_i^{(k)}(e^{-z/k})
 =
\exp\bigg\{-\sum_{i=1}^{[kq]}\frac{1}{k^{2}\gamma_k\zeta_i^{(k)}}
\psi_{\frac{i}{k}}^{(k)}(z)\bigg\},
 \eeqnn
where $\zeta_i^{(k)}\in[h_i^{(k)}(e^{-z/k}),1]$. It is easy to show that
$\prod_{i=1}^{[kq]}h_i^{(k)}(e^{-z/k})$ converges to $1$ uniformly on
 $[0,l]$ for every $l\geq 0$ if Condition (4.B) holds, and hence for each
$1\leq i\leq [kq]$, $\prod_{j=i+1}^{[kq]}h_j^{(k)}(e^{-z/k})$ converges to
$1$ uniformly on $[0,l]$ for every $l\geq 0$. By letting $k\to\infty$ in
(\ref{3.5}) we see $\phi_q^{(k)}(z)$ uniformly converge to a function
$\phi_q (z)$ on {$[0,l]$ for every $l\geq 0$} and (\ref{1.3}) holds. Then the desired result follows readily.\qed

\bproposition\label{p3.2} To each admissible family of branching
mechanisms $\{\phi_q: q\geq 0\}$ with $(\partial/\partial\theta)
\phi_{\theta}(z)=-\psi_{\theta}(z)$, there correspond two sequences
$\{\phi_{0}^{(k)}\}$ and $\{\psi_{\theta}^{(k)}\}$ in form of (\ref{4.4a})
and (\ref{4.4b}), respectively, so that Conditions (4.A) and (4.B) are
satisfied. \eproposition

\proof By Li (2011, p.93) there is a sequence $\{\phi_{0}^{(k)}\}$ in form
of (\ref{4.4a}) satisfying Condition (4.A). By Li (2011, p.102), there is
a family of probability generating functions $\{\bar{h}_{\theta}^{(k)}\}$
such that
 \beqnn
k[1-\bar{h}_{\theta}^{(k)}(e^{-z/k})]\to \psi_{\theta}(z)
 \eeqnn
uniformly on $[0,l]^2$ for every $a\geq 0$ as $k\to \infty$. Let
 \beqnn
\tilde{h}_{\theta}^{(k)}(z)=1+\frac{1}{k\gamma_k}[\bar{h}_{\theta}^{(k)}(z)-1],\quad
\theta\geq 0, \quad |z|\leq 1.
 \eeqnn
Clearly, $\{\tilde{h}_{\theta}^{(k)}:\theta\geq 0\}$ is a family of
probability generating functions and
 \beqnn
k^2\gamma_k [1-\tilde{h}_{\theta}^{(k)}(e^{-z/k})]\to \psi_{\theta}(z)
 \eeqnn
uniformly on $[0,l]^2$ for every $l\geq 0$ as $k\to \infty$. For each
$k\geq 1$, define $h_{i}^{(k)}=\tilde{h}_{i/k}^{(k)}$, $i=1,2,\cdots$.
Then by the continuity of $(\theta,z)\mapsto\psi_{\theta}(z)$ we get the
result.\qed

Given a constant $a\geq 0$, denote by $\{Y_{t,a}^{(k)}:t\geq 0\}$ the
restriction of $\{Y_t^{(k)}:t\geq 0\}$ to $[0,a]$.
 Then it is easy to see
 \beqnn
Y_{0,a}^{(k)}=\frac{1}{k}\sum_{i=0}^{[ka]}X_0^{(k)}(i)\delta_{\frac{i}{k}}
\quad\text{and}\quad
Y_{t,a}^{(k)}=\frac{1}{k}\sum_{i=0}^{[ka]}X_{[\gamma_kt]}^{(k)}(i)\delta_{\frac{i}{k}}.
 \eeqnn
Then $\{Y_{t,a}^{(k)}:t\geq 0\}$ is a measure-valued Markov process with
state space $M_k[0,a]$. From (\ref{3.3}) one can see the (discrete)
generator $L_k$ of $\{Y_{t,a}^{(k)}:t\geq 0\}$ is given by, for
$\nu=k^{-1}\sum_{i=0}^{[ka]}x_i^{(k)}\delta_{i/k}\in M_k[0,a]$ and $f\in
C[0,a]^{++}$,
 \beqlb\label{3.6}
L_ke^{-\<\nu,f\>} \ar=\ar {\gamma_k}\Big[\prod_{i=0}^{[ka]}
g_i^{(k)}(e^{-f(\frac{i}{k})/k})^{x_i}h_i^{(k)}(e^{-f(\frac{i}{k})/k})^{\sum_{j=0}^{i-1}x_j}
-e^{-\<\nu,f\>}\Big]\cr \ar=\ar
e^{-\<\nu,f\>}{\gamma_k}\Big[\exp\Big\{\sum_{i=0}^{[ka]}
\log\Big(g_i^{(k)}(e^{-f(\frac{i}{k})/k})^{x_i}h_i^{(k)}(e^{-f(\frac{i}{k})/k})^{\sum_{j=0}^{i-1}x_j}\Big)\cr
\ar\ar \qqquad\qqquad + \<\nu,f\>\Big\}-1\Big]\cr \ar=\ar
e^{-\<\nu,f\>}{\gamma_k}\Big[\exp\{\alpha_k+\beta_k\} -1\Big],
 \eeqlb
where
 \beqnn
\alpha_k =\sum_{i=0}^{[ka]}x_i\Big[\log
g_i^{(k)}(e^{-f(\frac{i}{k})/k})+f(\frac{i}{k})/k\Big], \quad\beta_k=
\sum_{i=0}^{[ka]}\sum_{j=0}^{i-1}x_j\log h_i^{(k)}(e^{-f(\frac{i}{k})/k}).
 \eeqnn
By the definition of $g_i^{(k)}$ we have
 \beqnn
\alpha_k \ar=\ar \sum_{i=0}^{[ka]}x_i\Big[\log
g_0^{(k)}(e^{-f(\frac{i}{k})/k}) +\sum_{j=0}^{i}\log
h_j^{(k)}(e^{-f(\frac{i}{k})/k})+f(\frac{i}{k})/k\Big]\cr \ar=\ar
\sum_{i=0}^{[ka]}x_i\Big[\log
g_0^{(k)}(e^{-f(\frac{i}{k})/k})+f(\frac{i}{k})/k\Big]
+\sum_{i=0}^{[ka]}\sum_{j=0}^{i}x_i\log h_j^{(k)}(e^{-f(\frac{i}{k})/k}).
 \eeqnn
It follows that
 \beqnn
\alpha_k+\beta_k \ar=\ar \sum_{i=0}^{[ka]}x_i\Big[\log
g_0^{(k)}(e^{-f(\frac{i}{k})/k})+f(\frac{i}{k})/k\Big]
+\sum_{i=0}^{[ka]}\sum_{j=0}^{i}x_i\log
h_j^{(k)}(e^{-f(\frac{i}{k})/k})\cr \ar\ar
+\sum_{i=0}^{[ka]}\sum_{j=0}^{i-1}x_j\log
h_i^{(k)}(e^{-f(\frac{i}{k})/k})\cr \ar=\ar \sum_{i=0}^{[ka]}x_i\Big[\log
g_0^{(k)}(e^{-f(\frac{i}{k})/k})+f(\frac{i}{k})/k\Big]
+\sum_{i=0}^{[ka]}\sum_{j=0}^{i}x_i\log
h_j^{(k)}(e^{-f(\frac{i}{k})/k})\cr \ar\ar
+\sum_{i=0}^{[ka]}\sum_{j=i+1}^{[ka]-1}x_i\log
h_j^{(k)}(e^{-f(\frac{j}{k})/k})\cr \ar=\ar \sum_{i=0}^{[ka]}x_i\Big[\log
g_0^{(k)}(e^{-f(\frac{i}{k})/k})+f(\frac{i}{k})/k\Big]
+\sum_{i=0}^{[ka]}\sum_{j=0}^{[ka]}x_i\log h_j^{(k)}(e^{-f(\frac{i\vee
j}{k})/k})\cr \ar=\ar \frac{1}{\gamma_k}\bigg[
\sum_{i=0}^{[ka]}\frac{x_i}{k\zeta_i^{(k)}}\phi_0^{(k)}(f(\frac{i}{k}))
-\sum_{i=0}^{[ka]}\frac{x_i}{k}\Big(\sum_{j=0}^{[ka]}\frac{1}{k\zeta_{i,j}^{(k)}}
\psi_{\frac{j}{k}}^{(k)}({f(\frac{i\vee j}{k}}))\Big)\bigg],
 \eeqnn
where $\zeta_i^{(k)}$ is between $e^{-f(\frac{i}{k})/k}$ and
$g_0^{(k)}(e^{-f(\frac{i}{k})/k})$, ~$\zeta_{i,j}^{(k)}\in
[h_j^{(k)}(e^{-f(\frac{i\vee j}{k})/k}),1]$. Clearly, both $\zeta_i^{(k)}$
and $\zeta_{i,j}^{(k)}$ converge to $1$ uniformly as $k\to\infty$ if
Conditions (4.A) and (4.B) hold. Then the above equality implies
 \beqlb\label{3.7}
\alpha_k+\beta_k=\frac{1}{\gamma_k}\Big[\<\nu,\phi_0^{(k)}(f(\cdot))\>-
\<\nu,\Psi^{(k)}(\cdot,f)\>+o(1)\Big],
 \eeqlb
where
 \beqnn
\Psi^{(k)}(\cdot,f)
=\sum_{j=0}^{[ka]}\frac{1}{k}\psi_{\frac{j}{k}}^{(k)}({f(\cdot\vee\frac{j}{k}})).
 \eeqnn

Let $\{Y_{t,a}:t\geq 0\}$ be the c\`{a}dl\`{a}g superprocess with
transition semigroup $(Q_t)_{t\geq 0}$ defined by (\ref{1.5}) and
(\ref{1.6}).
\btheorem\label{t3.2} Suppose that Conditions (4.A) and (4.B) are
satisfied.  If $Y_{0,a}^{(k)}$ converges to $Y_{0,a}$ in
distribution on $M[0,a]$, then $\{Y_{t,a}^{(k)}:t\geq 0\}$ converges to
$\{Y_{t,a}:t\geq 0\}$ in distribution on $D([0,\infty),M[0,a])$.\etheorem

\proof As in the proof of Theorem~2.1 in Li (2006), we shall prove the
convergence of the generators. Let $D_1$ be the algebra as defined in
Theorem \ref{t2.2}. For $f\in C[0,a]^{++}$, let
 \beqnn
Le^{-\<\nu,f\>}=e^{-\<\nu,f\>}\Big[\<\nu, \phi_0(f(\cdot))\>-\<\nu,
\Psi(\cdot,f)\>\Big], \qquad \nu\in M[0,a],
 \eeqnn
and extend the definition of $L$ to $D_1$ by linearity. By (\ref{1.5}) one
can check that $L$ is a restriction of strong generator of $(Q_t)_{t\geq
0}$. Note also that $L:=\{(f,Lf): f\in D_1\}$ is a linear space of
$C_0(M[0,a])\times C_0(M[0,a])$. On the other hand, let $f(x)=\lambda$ in
(\ref{1.5}) and (\ref{1.6}), we have the function $\lambda\mapsto
V_t(\lambda)$ is strictly increasing on $[0,\infty)$ for every $t\geq 0$;
see, e.g., Li (2011, p.58). Therefore, $V_t(\lambda)>0$ for every $\lambda
>0$ and $t\geq 0$. In view of (\ref{1.5}), for any $f\in C[0,a]^{++}$ we
have $V_tf\in C[0,a]^{++}$ for every $t\geq 0$. Then $D_1$ is invariant
under $(Q_t)_{t\geq 0}$,
which is a core of the strong generator of $(Q_t)_{t\geq 0}$;
see, e.g., Ethier and Kurtz (1986, p.17). In other words, the closure of
$L$ generates $(Q_t)_{t\geq 0}$ uniquely; see, e.g., Ethier and Kurtz
(1986, p.15 and p.17). Based on (\ref{3.6}) and (\ref{3.7}) one can see
 \beqnn
\lim_{k\to\infty}\sup_{\nu\in
M_k[0,a]}|L_ke^{-\<\nu,f\>}-Le^{-\<\nu,f\>}|=0
 \eeqnn
for every $f\in C[0,a]^{++}$,  which implies
 \beqnn
\lim_{k\to\infty}\sup_{\nu\in M_k[0,a]}|L_kF(\nu)-LF(\nu)|=0
 \eeqnn
for every $F\in D_1$. By Ethier and Kurtz (1986, p.226 and pp.233-234) we
conclude that $\{Y_{t,a}^{(k)}:t\geq 0\}$ converges to the immigration
superprocess $\{Y_{t,a}:t\geq 0\}$ in distribution on
$D([0,\infty),M[0,a])$. \qed

 Let
$\{0\leq a_1<a_2<\cdots<a_n=a\}$ be an ordered set of constants. Denote by
$\{Y_{t,a_i}:t\geq 0\}$ and $\{Y_{t,a_i}^{(k)}:t\geq 0\}$ the restriction of $\{Y_t:t\geq 0\}$ and $\{Y_t^{(k)}:t\geq 0\}$ to $[0,a_i]$, respectively.  Let
$Y_{t}(a_i):=Y_t[0,a_i]$ and $Y_{t}^{(k)}(a_i):=Y_t^{(k)}[0,a_i] $ for every $t\geq 0$, $i=1,2,\cdots,n$. By arguments similar to those in Section 3 we have:

\btheorem\label{t3.3} Suppose that Conditions (4.A) and (4.B) are
satisfied.  If $Y_{0,a}^{(k)}$
converges to $Y_{0,a}$ in distribution on $M[0,a]$, then
$\{(Y_{t,a_1}^{(k)}, \cdots, Y_{t,a_n}^{(k)}):t\geq 0\}$ converges to
$\{(Y_{t,a_1}, \cdots, Y_{t,a_n}):t\geq 0\}$ in distribution on
$D([0,\infty),M[0,a_1]\times \cdots\times M[0,a_n])$.\etheorem

\bcorollary\label{c3.4} Suppose that Conditions (4.A) and (4.B) are
satisfied. If  $(Y_{0}^{(k)}(a_1), \cdots, Y_{0}^{(k)}(a_n))$
converges to $(Y_{0}(a_1), \cdots, Y_{0}(a_n))$ in distribution on
$\mbb{R}_+^{n}$, then $\{(Y_{t}^{(k)}(a_1), \cdots,
Y_{t}^{(k)}(a_n)):t\geq 0\}$ converges to $\{(Y_{t}(a_1), \cdots,
Y_{t}(a_n)):t\geq 0\}$ in distribution on $D([0,\infty), \mbb{R}_+^{n})$.
\ecorollary

\bigskip

\textbf{Acknowledgements.} This paper was written under the supervision of
Professor Zenghu Li, to whom gratitude is expressed. We also would like to
acknowledge the Laboratory of Mathematics and Complex Systems (Ministry of
Education, China) for providing us the research facilities.

\bigskip\bigskip


\bigskip\bigskip

\noindent School of Mathematical Sciences

\noindent Beijing Normal University

\noindent Beijing 100875, P.\,R.\ China

\noindent E-mails: {\tt hehui@bnu.edu.cn, marugang@mail.bnu.edu.cn}

\end{document}